\PassOptionsToPackage{hyphens}{url}
\documentclass[graybox]{svmult}

\usepackage{type1cm}        
\usepackage{makeidx}         
\usepackage{graphicx}        
\usepackage{multicol}        
\usepackage[bottom]{footmisc}
\usepackage{epstopdf,psfrag}

\usepackage{newtxtext}       %
\usepackage[varvw]{newtxmath}       

\usepackage{graphicx}       

\usepackage[latin1]{inputenc}
\usepackage[T1]{fontenc}
\usepackage{microtype} 

\usepackage{cite}
\usepackage{array,colortbl}
\usepackage{amsmath,amsfonts,bm,mathtools,mathrsfs} 

\DeclareFontFamily{U}{mathx}{\hyphenchar\font45}
\DeclareFontShape{U}{mathx}{m}{n}{
      <5> <6> <7> <8> <9> <10>
      <10.95> <12> <14.4> <17.28> <20.74> <24.88>
      mathx10
      }{}
\DeclareSymbolFont{mathx}{U}{mathx}{m}{n}
\DeclareFontSubstitution{U}{mathx}{m}{n}
\DeclareMathAccent{\widecheck}{0}{mathx}{"71}

\def\citep#1#2{\cite[{#1}]{#2}}

\usepackage{algorithm}
\usepackage{algorithmic}

\newcommand{\bbE}{{\mathbb{E}}}

\newcommand{\R}{{\mathbb{R}}} 
\newcommand{\EE}{{\mathbb{E}}}

\newcommand{\bst}{{\boldsymbol{t}}}

\newcommand{\bsx}{{\boldsymbol{x}}}

\newcommand{\bsT}{{\boldsymbol{T}}}

\newcommand{\bsX}{{\boldsymbol{X}}}

\newcommand{\bszero}{{\boldsymbol{0}}} 

\newcommand{\bsPsi}{{\boldsymbol{\Psi}}}

\newcommand{\calO}{{\mathcal{O}}}

\newcommand{\calT}{{\mathcal{T}}}
\newcommand{\calU}{{\mathcal{U}}}

\newcommand{\abs}[1]{\left\vert#1\right\vert}

\usepackage{cite}
\usepackage{xspace}
\usepackage[bookmarks=false]{hyperref}
\usepackage[hyphens]{url}

\begin{document}

\title*{Challenges in Developing Great Quasi-Monte Carlo Software}
\authorrunning{S.-C.\ T.\ Choi et al.}
\author{Sou-Cheng T.~Choi \and Yuhan Ding \and Fred J. Hickernell \and Jagadeeswaran  Rathinavel \and Aleksei G. Sorokin}
\institute{
Sou-Cheng T.~Choi \at Department of Applied Mathematics, Illinois Institute of Technology; SAS Institute Inc. \\RE 220, 10 W.\ 32$^{\text{nd}}$ St., Chicago, IL 60616 \email{schoi32@iit.edu}
\and
Yuhan Ding \at
 Department of Applied Mathematics, Illinois Institute of Technology,\\ RE 220, 10 W.\ 32$^{\text{nd}}$ St., Chicago, IL 60616 \email{yding2@hawk.iit.edu}
\and
Fred J. Hickernell \at Center for Interdisciplinary Scientific Computation and \\
Department of Applied Mathematics, Illinois Institute of Technology \\ RE 220, 10 W.\ 32$^{\text{nd}}$ St., Chicago, IL 60616 \email{hickernell@iit.edu}
\and
Jagadeeswaran  Rathinavel \at
Department of Applied Mathematics, Illinois Institute of Technology,\\ RE 220, 10 W.\ 32$^{\text{nd}}$ St., Chicago, IL 60616 \email{jrathin1@hawk.iit.edu}
\and
Aleksei G. Sorokin \at
Department of Applied Mathematics, Illinois Institute of Technology,\\ RE 220, 10 W.\ 32$^{\text{nd}}$ St., Chicago, IL 60616 \email{asorokin@hawk.iit.edu}}

\maketitle

\abstract{
Quasi-Monte Carlo (QMC) methods have developed over several decades. With the explosion in computational science, there is a need for great software that implements QMC algorithms. We summarize the  QMC software that has been developed to date, propose some criteria for developing great QMC software, and suggest some steps toward achieving great software.  We illustrate these criteria and steps with the Quasi-Monte Carlo Python library (QMCPy), an open-source community software framework, extensible by design with common programming interfaces to an increasing number of existing or emerging QMC libraries developed by the greater community of QMC researchers.
}

\section{Introduction}

A QMC approximation of $\mu := \EE[f(\bsX)]$, $\bsX \sim \calU[0,1]^d$ can be implemented in a few steps:
\begin{enumerate}
    \item Draw a sequence of $n$ low discrepancy (LD) \cite{Nie92, SloJoe94, DicEtal14a, DicEtal22a} nodes,
    $\bsx_1,\dots,\bsx_n$ that mimic $\calU[0,1]^d$.
    \item Evaluate the integrand $f$ at these nodes to obtain $f(\bsx_i)$, $i=1,\dots,n$.
    \item Estimate the true mean, $\mu$, by the sample mean,
    \begin{equation*}
        \hat{\mu}_n := \frac{1}{n} \sum_{i=1}^n f(\bsx_i).
    \end{equation*}
\end{enumerate}

However, the practice of QMC is often more complicated.  The original problem may need to be rewritten to fit the above form and/or to facilitate a good approximation with a small number of sampling points, $n$.  Practitioners may wish to determine $n$ adaptively to satisfy a prescribed error tolerance.

In the next section, we describe why QMC software is important.  We then discuss the characteristics of great QMC software:

\begin{itemize}
\item Integrated with related libraries (Sect.\ \ref{CDHJS_sec:integrated}),

\item Correct (Sect.\ \ref{CDHJS_sec:correct}),

\item Efficient in computational time and memory (Sect.\ \ref{CDHJS_sec:efficient}),

\item Accessible to practitioners and theorists alike (Sect.\ \ref{CDHJS_sec:accessible}), and

\item Sustainable by a community that owns it (Sect.\ \ref{CDHJS_sec:sustainable}).

\end{itemize}
In each section, we describe the challenges faced and how they might be overcome.

We draw on our collective experiences as
members of the academic QMC community, developers of open-source and commercial scientific software, and users of a wide range of scientific software as research scientists or data scientists. Many of the insights that we have gained have come through our development of the Guaranteed Automatic Integration Library (GAIL)~\cite{TonEtAl22a} in MATLAB and the Quasi-Monte Carlo Python library (QMCPy)~\cite{QMCPy2020a}.

\section{Why Develop QMC Software} \label{CDHJS_sec:why_we_need_software}

Recent interest in quality scientific software is exemplified by the 2021 US Department of Energy report, \emph{Workshop on the Science of Scientific-Software Development and Use} \cite{ASCR-SSSDU, osti_1846008}. This workshop not only discussed diagnostics and treatments for the challenges of developing great scientific software, but also emphasized the importance of implementing theoretical advancements into well-written, accessible software libraries. Three cross-cutting themes arose in this workshop:
\begin{itemize}
    \item  We need to consider both human and technical elements to better understand how to improve the development and use of scientific software.

    \item We need to address urgent challenges in workforce recruitment and retention in the computing sciences with growth through expanded diversity, stable career paths, and the creation of a community and culture that attract and retain new generations of scientists.

    \item Scientific software has become essential to all areas of science and technology, creating opportunities for expanded partnerships, collaboration, and impact.
\end{itemize}
These themes apply to QMC software in particular, as well as scientific software in general.

\subsection{QMC Theory to Software}

QMC software makes theoretical advances in QMC available to practitioners. However, translating pseudo-code into good executable code is nuanced.  Software developers must write code that is computationally efficient, numerically stable, and provides reasonable default choices of tuning parameters. Great QMC software relieves users of these concerns.

Great QMC software eliminates the need for researchers to resort to  ``do-it-yourself'' for established algorithms.  The less code we have to write, the fewer the errors.

\subsection{QMC Software to Theory}

Great QMC software opens up new application areas for QMC methods by allowing practitioners to compare new methods to existing ones.  QMC software has been successful in quantitative finance, uncertainty quantification, and image rendering. Unexpectedly good or bad computational results prompt theoretical questions.  For example, the early application of QMC to high dimensional integrals arising in financial risk \cite{PasTra95} led to a wave of theoretical results on the tractability of integration in weighted spaces \cite{Woz99a,DicEtal14a,NovWoz10a}.

\bigskip

We next discuss the characteristics of great QMC software as outlined in the introduction.  For each characteristic, we identify potential deficiencies and suggest possible remedies.

\section{Integrated} \label{CDHJS_sec:integrated}

Expanding on the problem formulation in the introduction, we want to approximate well the integral or expectation, $\mu$, by the sample mean, $\hat{\mu}_n$:
\begin{gather}
\label{CDHJSeq:cubSummary}
	\mu : = \int_\calT g(\bst) \, \lambda(\bst) \, \D \bst  = \bbE[f(\bsX)] = \int_{[0,1]^d} f(\bsx)  \, \D \bsx \approx \frac 1n \sum_{i=1}^{n} f(\bsx_i) =: \hat{\mu}_n, \\
 \label{CDHJS_eq:error_tol}
 \abs{\mu - \hat{\mu}_n} \le \varepsilon.
\end{gather}
Progressing from the original problem to a satisfactory solution requires several software components, which we describe in the subsections below.  Great QMC software environments allow different implementations of these components to be freely interchanged.

\subsection{LD Sequence Generators} 
\begin{figure}[t]
    \centering
    \includegraphics[width=\textwidth]{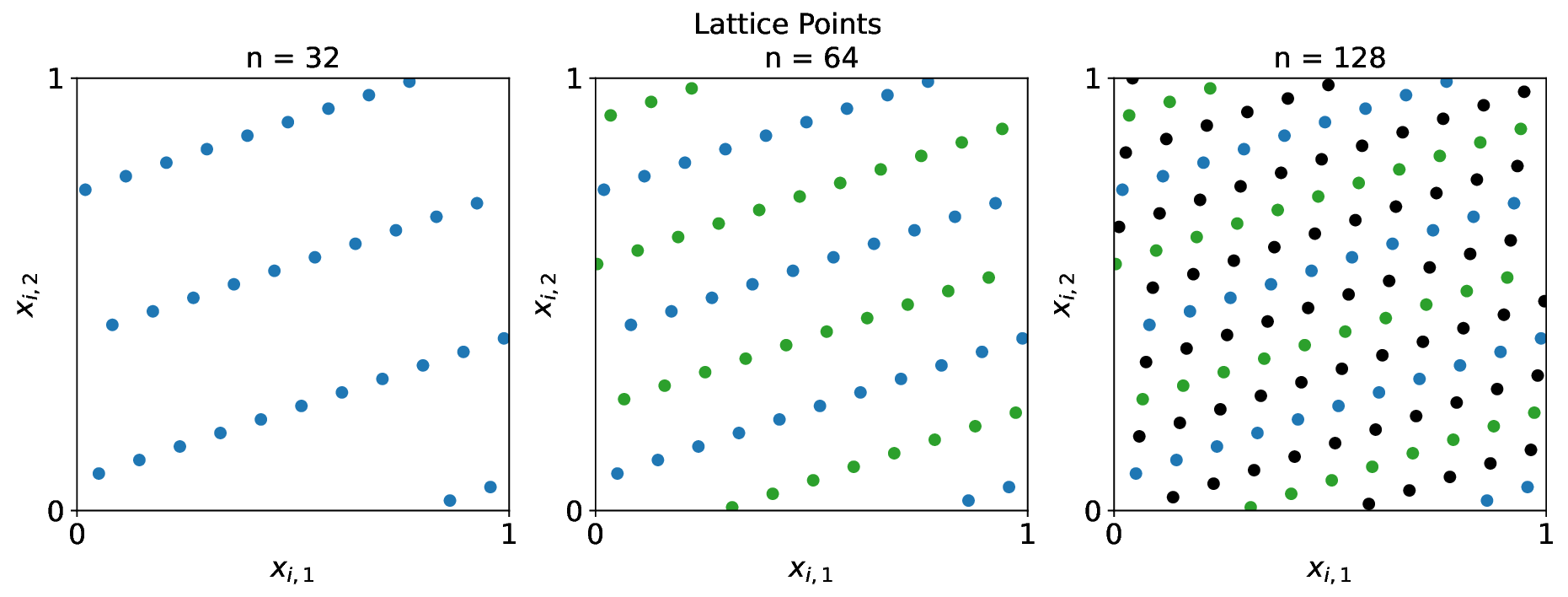}
    \caption{A two-dimensional projection of lattice points for increasing sample size.  The later points fill in gaps because they are highly correlated with the earlier points.  These plots may be reproduced using the Jupyter notebook \cite{MCQMC2022Figs}. }
    \label{CDHJSfiglatticepts}
\end{figure}

A lattice, digital sequence, or Halton sequence generator typically supplies the sequence $\bsx_1, \bsx_2, \ldots$.  These sequences may be deterministic or random and are intended to have an empirical distribution that approximates well the uniform distribution.  Figure \ref{CDHJSfiglatticepts} illustrates an example of a (randomly shifted) lattice sequence.

LD generators are the most prevalent, generally available components of QMC software.  These generators appear in Association of Computing Machinery  publications \cite{BraFox88,BraFoxNie92,HonHic00a}, FinDer \cite{PasTra95,FinDer}, libseq \cite{FriKel02,FriKelweb}, BRODA \cite{BRODA20a}, NAG \cite{NAG27}, MATLAB \cite{MAT9.13}, SamplePack \cite{SamplePack}, Gr\"unschlo{\ss}'s website \cite{GruWeb}, the Magic Point Shop \cite{NuyWeb}, R \cite{QRNG2020}, Julia \cite{QMCJulia}, PyTorch \cite{paszke2019pytorch}, SAS \cite{SAS_LD}, SciPy \cite{virtanen2020scipy},  TensorFlow \cite{tfqfQMC2021a}, MatBuilder \cite{paulin2022}, and  QMCPy \cite{QMCPy2020a}.

Lattice and digital sequences are not unique, and improved or alternate versions continue to be constructed.  Stephen Joe and Frances Kuo \cite{JoeKuo03,KuoJoe08a,SobolDirection} have proposed  Sobol' generator direction numbers, which are now widely used, but were not part of early LD sequence software.  Pierre L'Ecuyer and collaborators have created LatNet Builder \cite{LEcEtal22a} to construct good lattice and polynomial lattice sequences based on user-defined criteria. During lunch at MCQMC 2022, it was proposed that we  agree on a consistent format for storing the parameters that define good LD sequences so that they can be shared across software libraries.

Users computing solutions in comprehensive software environments, such as NAG, MATLAB, R, SciPy, and TensorFlow, have convenient access to LD sequence generators.  However, most of these implementations in large libraries would benefit a wider menu of LD offerings, including custom generators from LatNet Builder or other sources.  Some large environments such as  Dakota \cite{DakotaUsersManual} have quite limited LD generators and should expand their offerings. SAS facilitates flexible integration of third-party open-source software; see \href{https://www.sas.com/en_us/software/viya/open.html}{\url{https://www.sas.com/en_us/software/viya/open.html}}.
This presents opportunities for developers of smaller libraries with more QMC features, like BRODA and QMCPy, to connect their libraries with larger software environments. An example of QMCPy connecting with other software is given in the next subsection.

In fact, QMCPy is designed not simply to contain our own algorithms.  It is designed to serve as a middleware framework, allowing users to tap into a collection of QMC libraries originally written in various languages by different researchers at different times. Using the object-oriented design, QMCPy specifies a few  extensible high-level mathematical objects as described in this section and in \cite{ChoEtal22a}. Developers of other (Q)MC software can  extend these objects to incorporate their existing algorithms by a combination of mechanisms: object inheritance, Python package inclusion, direct code translation, and/or interfaces to compiled C or Java libraries.

\subsection{Integrands and Variable Transformations} \label{CDHJS:sec:integrands}
The original integral or expectation arising from an application is defined in terms of the integrand, $g :\mathcal{T} \to \R$, and the non-negative weight, $\lambda$, as described in \eqref{CDHJSeq:cubSummary}.  For example,
\begin{itemize}
	\item $g$ is the discounted payoff of a financial derivative and $\lambda$ is the probability density function (PDF) for a discretized Brownian motion \cite{Gla03},

	\item  $g$ is the velocity of a fluid at a point in rock with a random porosity field whose discretized PDF is $\lambda$ \cite{KuoNuy16a}, or

	\item $g\lambda$ is the unnormalized Bayesian posterior density multiplied by the parameter of interest \cite{GelEtal13}.
\end{itemize}

To put $\mu$ into the form that is directly accessible to QMC methods requires a variable transformation, $\bsPsi$, such that
\begin{equation*}
    \bsPsi((0,1)^d) = \calT, \qquad f(\bsx) = g(\bsPsi(\bsx)) \lambda(\bsPsi(\bsx)) \abs{\frac{\partial \bsPsi}{\partial\bsx}}.
\end{equation*}
Such a transformation is typically non-unique, and the choice of a good one is equivalent to importance sampling.  If $\lambda$ is the PDF for a random variable $\bsT$, then $\mu = \EE[g(\bsT)]$.  In this case, one might choose $\bsPsi$ such that $\bsT \sim \bsPsi(\bsX)$, $\bsX \sim \calU(0,1)^d$, in which case $\lambda(\bsPsi(\bsx)) \abs{\partial \bsPsi/\partial\bsx} = 1$ and $f(\bsx) = g(\bsPsi(\bsx))$. Software such as SciPy and QMCPy facilitate such variable transformations by providing default $\bsPsi$ for common $\bsT$. This alleviates the burden on the user to determine a transformation and implement the Jacobian.

\begin{figure}[t]
    \centering
    \raisebox{0.5cm}{\includegraphics[width=0.4\textwidth]{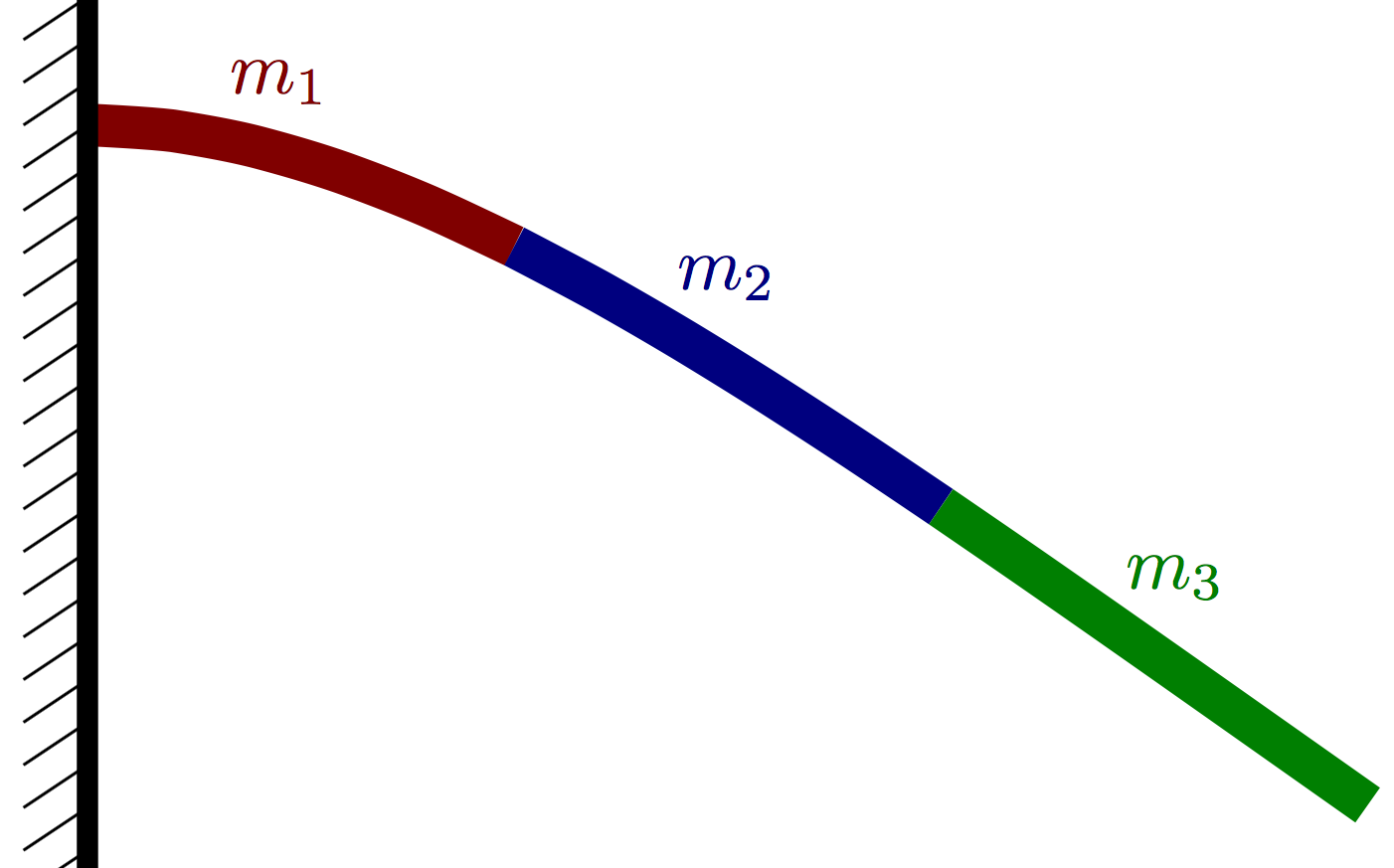}}
    \includegraphics[width=0.5\textwidth]{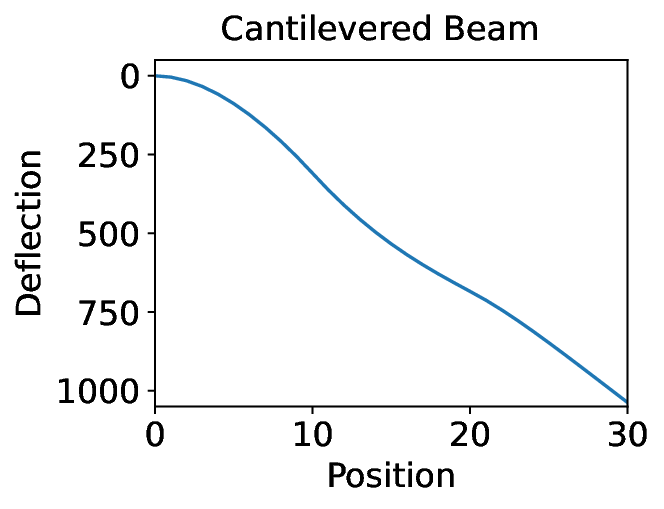}
    \includegraphics[width=1\textwidth]{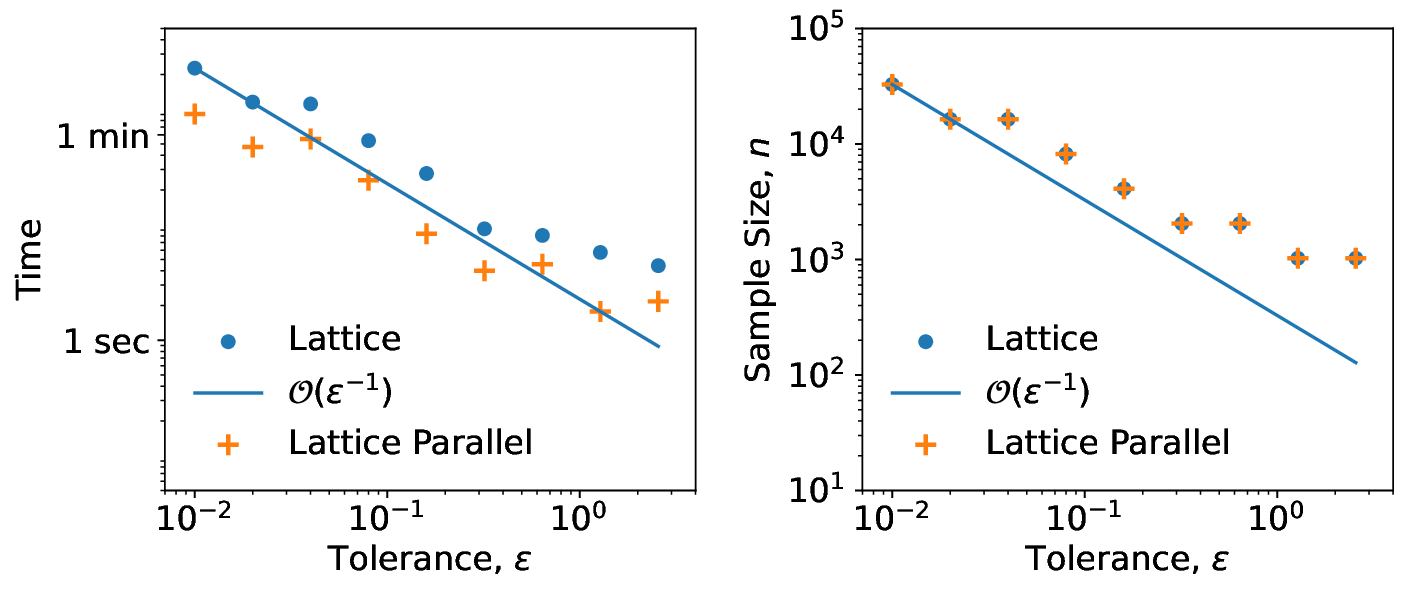}
    \caption{The upper left plot shows a piecewise log effective stiffness (image reproduced from~\cite{ParSee22a} with permission). The upper right plot shows the expected deflection of a cantilevered beam with random log effective stiffness. In the lower plots, QMCPy computes the answer in serial and in parallel using UM-Bridge and docker containers.  For small error tolerances, the time and sample size are $\calO(\varepsilon^{-1})$ and parallel processing gives a speed-up factor of four to five. These plots may be reproduced using the Jupyter notebook \cite{MCQMC2022Figs}.  }
    \label{CDHJSfigbeam}
\end{figure}

Complex integrands might be best evaluated by libraries that have no direct connection to QMC software libraries.  In the uncertainty quantification (UQ) space, the recently developed UQ and Model Bridge (UM-Bridge) \cite{UMBridge23}, connects libraries that generate sequences of model parameters, $\bsx_1, \bsx_2, \ldots$ with libraries that evaluate models, $f(\bsx_1), f(\bsx_2), \ldots$, and their derivatives if applicable. UM-Bridge allows coupling of UQ codes to models regardless of the programming languages.

We illustrate in Figure~\ref{CDHJSfigbeam} how UM-Bridge connects QMCPy with a differential equation boundary value problem solver for a cantilevered beam of length~$L$, with a random log effective stiffness $m(\cdot)$. Following ~\cite{ParSee22a}, we assume $m(\cdot)$ is piecewise constant over three regions and the three constants take a $\calU[1,1.2]^3$ distribution.
The beam position $x$ is discretized at $N=31$ equidistant points, and
the vertical deflection $u(x)$ when the beam is subject to a given external load  $g(x)$ satisfies the Euler-Bernoulli partial differential equation (PDE):
$$ \frac{\partial^2}{\partial x^2}\left[ \exp(m(x)) \frac{\partial^2  u(x)}{\partial x^2}\right] = g(x),$$
with the following boundary conditions:
$$ \left.u(x)\right|_{x=0} = 0, \qquad
\left.\frac{\partial u(x)}{\partial x}\right|_{x=0} = 0, \qquad
\left.\frac{\partial^2 u(x)}{\partial x^2}\right|_{x=L} = 0, \qquad  \left.\frac{\partial^3 u(x)}{\partial x^3}\right|_{x=L} = 0. $$
This PDE can be discretized with finite differences yielding a linear system $K(m) \hat{u} = \hat{g}$, where $\hat{u}$ and $\hat{g} \in \mathbb{R}^N$.
The beam model was written in C++ using the MUQ's model framework~\cite{MUQ2}. 

QMCPy simulates the random Young's modulus and computes the expected deflection of the beam at the node points for a sequence of error tolerances using lattice rules.  By taking advantage of clusters and containers, UM-Bridge provides automatic load balancing to facilitate seamless parallel model evaluation. As a result, thread-parallel QMCPy can offload costly model evaluations to remote clusters~\cite{seelinger2023lowering}.

We need more connections of QMC software to other libraries, like what has been done with  QMCPy and UM-Bridge.  This will allow newer and better QMC algorithms to tackle new applications.

\subsection{Stopping Criteria}
The cantilevered beam example illustrated in Figure \ref{CDHJSfigbeam} utilizes QMCPy's stopping criterion to choose the sample size, $n$, based on the input error tolerance, $\varepsilon$, and the observed discrete Fourier series coefficients of the integrand, $f$.  See \cite{HicJim16a} for details of this stopping criterion for lattice rules and \cite{HicEtal18a} for a review on automatic QMC stopping criteria.  Theoretically justified stopping criteria for choosing sample sizes are important and a relatively recent development \cite{HicEtal14a,HicJim16a,JimHic16a,RatHic19a,JagHic22a}.  Some of these stopping criteria for absolute error criteria, \eqref{CDHJS_eq:error_tol},  have been extended to relative error criteria or a combination of absolute and relative error criteria, as well as stopping criteria for approximating functions of several integrals \cite{Jia16a,GilJim16b,HicEtal17a,JagSor23a}.

\bigskip

We emphasize again that no one software library contains all of the best LD generators, integrands, variable transformations, and stopping criteria, which is why great QMC software must be well integrated into other software libraries.  This may take the form of large libraries incorporating more QMC routines and/or building connections among QMC libraries and with other non-QMC libraries.

\section{Correct} \label{CDHJS_sec:correct}
Although correctness is an obvious feature of great QMC software, it should not always be assumed in popular packages.  Here are two examples.

MATLAB's original randomization of Sobol' points was incorrectly implemented.  This was discovered by Llu\'is Antoni Jim\'enez Rugama, who informed the developer.  The error was corrected in MATLAB R2018a.

SciPy and PyTorch's original implementations of Sobol' sequences omitted the first point.  The rationale was that the unscrambled first point of the Sobol' sequence is $\bsx_1 = \bszero$, which becomes the undesirable $\boldsymbol{-\infty}$ when transformed to mimic a Gaussian distribution.  A vigorous discussion ensued \cite{scipySobol2020a,pytorchSoboldiscussion2020a}, which prompted an article by Art Owen \cite{Owe22a} on why one should not drop the first point and destroy the digital net structure.  Instead one should randomly scramble the Sobol' sequence (and other LD sequences) so that $\bsx_1 \ne \bszero$.  Although SciPy and PyTorch were corrected, other software libraries with LD sequence generators still omit the first point.

QMC software is sometimes written by those not intimately familiar with QMC.  Those of us most familiar with good QMC practices need to pay attention to QMC algorithms in popular software libraries and inform the developers when their software contains bugs or deviates from best practices.

The above examples and discussion also underscore the importance of continuous testing in QMC software development. This involves creating a comprehensive test suite that covers a wide range of scenarios and input parameters. Additionally, QMC software should be benchmarked against known results to ensure that it is both accurate and efficient.

\section{Efficient} \label{CDHJS_sec:efficient}
Another obvious feature of great QMC software is computational efficiency.  This requires developers to identify bottlenecks in performance. Intensive computations may need to be rewritten in lower-level languages, like C. We have found looped randomization and point generation routines to greatly benefit from low level implementations. Code may also need to be refactored to avoid excessive calls to memory.  As new hardware architectures arise, QMC software should adapt to take advantage of them.

Efficiency gains in QMC software can also result from clever algorithm choices. Adaptive QMC algorithms as in \cite{JimHic16a, HicJim16a, HicEtal17a, RatHic19a,JagHic22a} can be used to dynamically increase the number of sampling points based on estimated error bounds. They can lead to significant efficiency gains by reducing the number of unnecessary function evaluations. The Bayesian stopping rules for lattice and Sobol sequences in \cite{RatHic19a,JagHic22a} choose covariance kernels that match the corresponding LD sequence.  This reduces the computations from what would be $\calO(n^3)$ for arbitrary covariance kernels to only $\calO(n \log n)$ for these well-chosen kernels. Consequently, computing the stopping criterion is roughly the same order of magnitude as computing the solution.

In spite of the desire for efficiency, great QMC software will be written in a variety of computer languages since different communities tend to favor different languages. Trade-offs will be made between efficiency and convenience.
As noted by a referee, QMCPy is not particularly fast at generating tens of thousands of replicates of randomized QMC sample means in comparison to SSJ, which is based on Java. Such cases may point to the need for a Python layer over a faster implementation.

\section{Accessible} \label{CDHJS_sec:accessible}
One attribute of great QMC software that often receives insufficient attention is its accessibility.  Developers may write code for their own purposes and move on.

Great software should be stored in a repository, where the latest version can be downloaded.  There should be documentation that guides the user through installation and explains the library's features.  Demos are crucial for showing how the software works and for providing a template for new users to write their own code.  Tutorials help newcomers to QMC understand its advantages. Blogs educate a wider audience and enlarge our QMC community.  There should be a place for users to report bugs and make feature requests.

Accessible QMC software has greater impact.  Accessibility also depends on the next characteristic of great QMC software.

\section{Sustainable} \label{CDHJS_sec:sustainable}
Great QMC software libraries should be awake, not dormant.  They should be under active development or merged into other libraries that are actively maintained.  Commercial software has an infrastructure to maintain their algorithms, but these algorithms need periodic refreshing. Open source software requires a community of active users and developers to keep the library up-to-date and bug-free.

The QMC community can help sustain great QMC software in several ways.  \begin{itemize}
	\item Let's build our research code and showcase our new algorithms using existing QMC software libraries or build and maintain new libraries, if needed.  This allows the next generation to more easily reproduce our results and build upon them.

	\item Many of us QMC researchers rely on transient team members (students and postdocs) to write the code that illustrates our new ideas.  Building our own QMC libraries or contributing to existing ones helps ensure that the code our team members write lasts beyond their involvement in our teams.

	\item It is relatively easy to post demos of our new ideas on our own web pages, and if possible in the software libraries that we are using.  QMCPy welcomes blogs, short articles introducing a broader audience to important ideas, illustrated by software.

	\item When we see bugs or missing features in QMC software, let's persistently request fixes and improvements.  The more robust our QMC software, the easier it is to open new application areas to QMC and grow our QMC community.

    \item Let's publish our QMC  software with journals such as ACM TOMS, JORS, etc. to obtain code review and
    feedback.

    \item When we find QMC scientific software useful in our research or applications, let's encourage and recognize such software by properly citing it~\cite{smith2016software}.

    \item Encourage QMC community members to learn about software development through interdisciplinary education, internships, or hackathons.

\end{itemize}

Encouraging students and postdocs to invest time in developing QMC software will require us as their mentors valuing software as a research output.  Great software should be published with a digital object identifier (DOI).

The term research software engineer (RSE) is becoming a more valued vocation.  According to the website for  Research Software Engineers International,
\href{https://researchsoftware.org}{\url{https://researchsoftware.org}}, \emph{Research Software Engineers are people who combine professional software expertise with an understanding of research.}

We recognize the crucial role of non-code contributors to the project~\cite{K23}. Besides coding, we also need documentation, data visualization, testing, user interface or experience design, web content development, community management, and release management. These non-code contributions help make the project more  accessible. Non-coding roles offer opportunities for individuals with less coding experiences to develop technical and non-technical skills over time.  Additionally, the emergence of no-code or low-code analytical software such as RapidMiner, KNIME, Orange, SAS Studio, or SAS Viya enable individuals less experienced with complex coding to learn and contribute in a user-friendly environment.

\section{Summary and Future Work}
QMC has experienced decades of exciting growth in theory since the late 1950s.  The earliest publicly available QMC software dates to the early 1990s and its growth has not matched that of theory.  As a QMC community, we have a responsibility to implement our great theory into great QMC software.

In this paper, we have focused predominantly on employing QMC sampling in estimating the expectations of functions. However, it can be extended to other important mathematical problems such as function approximation \cite{Kmmerer2015, PotSch21}, simulation of Markov chains \cite{l2018sorting,puchhammer2021variance},  density estimation \cite{LEcuyer2022b}.
The SSJ library~\cite{l2002ssj} is a comprehensive tool kit that includes a wide range of QMC point sets, sequences, and randomizations in the ``hups'' and ``mcqmctools'' packages. It will be of great interest to the QMCPy community to make connections to these and other well established software in order to become applicable to an increasingly wider scope of problems and applications of computational mathematics and statistics.

Although Python is popular and versatile, it can be slower than highly optimized implementations in languages like Java or C++.  Developers of QMC software in Python would do well to investigate interfaces to optimize code in other languages or consider other code acceleration techniques provided by libraries such as Numba, JAX, Joblib, or Ray.  With time, languages rise and fall in popularity, and those who wish to provide great QMC software may need to pivot to new software environments.

\section*{Acknowledgements}
The authors are grateful to the referee's constructive comments. We also want to thank Linus Seelinger for feedback and discussion. The authors acknowledge many students and collaborators who have provided insight in our journey to develop great QMC software.



\end{document}